\documentclass[a4paper]{amsart}
\usepackage[12pt]{extsizes}
\usepackage{amsfonts,amsmath,amsthm,amscd,amssymb,latexsym,mathrsfs,enumitem}
\usepackage[english]{babel}
\usepackage{geometry}
\theoremstyle{plain}

\newtheorem{theorem}{Theorem}

\newtheorem{definition}{Definition}
\newtheorem{remark}{Remark}

\newtheorem{thmm}{Theorem}

\numberwithin{equation}{section}

\begin{document}

\title[A characterization of the uniform convergence points set]{A characterization of the uniform convergence points set of some convergent sequence of functions}

\author{Olena Karlova${}^{1,2}$}

\address{${}^1$ Jan Kochanowski University in Kielce, Poland; ${}^2$  Chernivtsi National University, Ukraine}
\email{maslenizza.ua@gmail.com}


\begin{abstract}
  We characterize the uniform convergence points set of a pointwisely convergent sequence of real-valued functions defined on a perfectly normal space.  We prove that if $X$ is a perfectly normal space which  can be covered by a disjoint sequence of dense subsets and  $A\subseteq X$, then $A$   is the set of points of the uniform convergence for some convergent sequence $(f_n)_{n\in\omega}$ of functions $f_n:X\to \mathbb R$ if and only if $A$ is  $G_\delta$-set which contains all isolated points of $X$. This result generalizes a theorem of J\'{a}n Bors\'{i}k published in 2019.
\end{abstract}
\maketitle

\bigskip

\begin{center}
  {\it Dedicated to the memory of J\'{a}n Bors\'{i}k}
\end{center}

\section{Introduction}
Let $X$ be a topological space, $(Y,|\cdot-\cdot |)$ be a metric space; $B(a,r)$ and $B[a,r]$ be an open and a closed ball in $Y$ with a center $a\in Y$ and a radius $r>0$, respectively.
By $\partial A$ we denote a boundary of a set $A$.

Let $\mathscr F=(f_n)_{n\in\omega}$ be a sequence of functions $f_n:X\to Y$. We denote   $PC(\mathscr F)$ the set of all points $x\in X$ such that the sequence $(f_n(x))_{n\in\omega}$ is convergent in $Y$. Therefore, we define the limit function $f(x)$ by the rule $f(x)=\lim_{n\to\infty} f_n(x)$ for all $x\in PC(\mathscr F)$. Let us observe that the set $PC(\mathscr F)$ can be represented in the form
\begin{gather}
  PC(\mathscr F)=\bigcap_{k\in\omega}\bigcup_{n\in\omega}\bigcap_{m\in\omega} f_{n+m}^{-1}(B[f(x),\tfrac{1}{k+1}]).
\end{gather}
If every function $f_n$ is continuous, then the set $PC(\mathscr F)$ is  $F_{\sigma\delta}$ in $X$. Hans Hahn \cite{Hahn} and Wac{\l}aw Sierpi\'{n}ski \cite{Sierp} proved independently that  the converse proposition is true for metrizable $X$ and $Y=\mathbb R$, that is, for every $F_{\sigma\delta}$-subset $A$ of a metrizable space $X$ there exists a sequence $\mathscr F$ of real-valued continuous functions $f_n:X\to\mathbb R$ such that $A=PC(\mathscr F)$.

After appearance of this theorem many results were obtained in similar directions: other types of convergence and other classes of functions were considered (see, for instance, \cite{Borsik2019,DST,Holi,MMS,NW1,NW2,W1,W2,W3}). J\'{a}n Bors\'{i}k  studied in \cite{Borsik2019}, in particular, the uniform convergence points set of a (convergent pointwisely) sequence of functions.

\begin{definition}{\rm
  A sequence $(f_n)_{n\in\omega}$ of functions $f_n:X\to Y$ between a topological space $X$ and a metric space $(Y,|\cdot - \cdot |)$ is
   {\it uniformly Cauchy  at a point $x_0\in X$}, if for every $\varepsilon>0$ there exist a neighborhood $U$ of $x_0$ and a number $n_0\in\omega$
   such that $|f_n(x)-f_m(x)|<\varepsilon$ for all $n,m\ge n_0$ and $x\in U$.}
\end{definition}

   Let $UC(\mathscr F)$ be a set of all points with the uniform Cauchy property for a sequence $\mathscr F=(f_n)_{n\in\omega}$.
  It is easy to see that if $(f_n)_{n\in\omega}$ is convergent pointwisely on $X$ to a function $f:X\to Y$, then
  \begin{gather}\label{gath1}
  UC(\mathscr F)=\{x\in X:\forall \varepsilon>0\,\, \exists U\ni x_0 \,\,\exists n_0\,\, \forall n\ge n_0 \,\,|f_n(y)-f(y)|<\varepsilon\,\, \forall y\in U\}.
  \end{gather}
Moreover, in this case $UC(\mathscr F)$ is the set of all points of the uniform convergence of $\mathscr F$.

Bors\'{i}k proved the following result.

\begin{thmm}[Bors\'{i}k, \cite{Borsik2019}]
  Let $X$ be a metric space and $A\subseteq X$. Then $A=UC(\mathscr F)$ for some convergent sequence   $\mathscr F=(f_n)_{n\in\omega}$ of functions $f_n:X\to \mathbb R$ if and only if   $A$ is $G_\delta$ and contains all isolated points of $X$.
\end{thmm}

This short note is inspirited by the above mentioned paper of J\'{a}n Bors\'{i}k. We generalize his theorem on a wider class of topological spaces.

\begin{definition}{\rm
A topological space $X$ is \emph{$\omega$-resolvable} if there exists a partition $\{X_n:n\in\omega\}$ of $X$ by  dense subsets.}
 \end{definition}

 For crowded spaces (i.e.,  spaces without isolated points) the class of all $\omega$-resolvable spaces includes all  metrizable spaces, Hausdorff countably compact spaces, arcwise connected spaces, etc.~\cite{CGF}

The main result of our note is the following theorem.

\begin{theorem}
  Let $X$ be a perfectly normal $\omega$-resolvable space and $A\subseteq X$. The following conditions are equivalent:
   \begin{enumerate}
     \item[(i)] $A$ is a set of all points of uniform convergence for some convergent sequence $(f_n)_{n\in\omega}$ of functions $f_n:X\to \mathbb R$;

     \item[(ii)] $A$ is a $G_\delta$-set which contains all isolated points of $X$.
   \end{enumerate}
\end{theorem}

\section{Proof of Theorem 1}

The implication {\textbf{\it (i) $\Rightarrow$ (ii)}} follows immediately from the equality~(\ref{gath1}).

{\textbf{\it (ii) $\Rightarrow$ (i)}}. Let $(G_n)_{n\in\omega}$ be a sequence of open sets in $X$ such that  $A=\bigcap_{n\in\omega} G_n$,  $G_{n+1}\subseteq G_n$ for every $n\in\omega$ and let $G_0=X$.

Since $X$ is perfectly normal, for every $n\in\mathbb N$ there exist   continuous functions $\varphi_n,\psi_n:X\to [0,1]$ such that $\varphi_n^{-1}(0)=\overline{G_n}$ and $\psi_n^{-1}(0)=X\setminus G_n$. Then every function $\alpha_n:X\to [-1,1]$ defined by the formula $\alpha_n=\varphi_n-\psi_n$ has the following properties:
\begin{eqnarray*}
\alpha_n(x)<0& \phantom{a}&\forall x\in G_n,\\
\alpha_n(x)=0& \phantom{a}&\forall x\in\partial G_n,\\
\alpha_n(x)>0& \phantom{a}&\forall x\in X\setminus \overline{G_n}.
\end{eqnarray*}

We consider functions  $\beta_n:X\to [-1,1]$ defined by the rule
$$
\beta_n(x)=\max_{i\le n}\alpha_i(x)
$$
for all $x\in X$. Then we claim that
$$
\partial G_n=\beta_n^{-1}(0)
$$
for every $n\in\mathbb N$. We need to prove   $\beta_n(x)=0 \,\,\Leftrightarrow\,\,\alpha_n(x)=0$ for every $x\in X$. Assume that $\beta_n(x)=0$. Then $\alpha_i(x)\le 0$ for all $i\le n$ and  $\alpha_k(x)=0$ for some $k\le n$. Then $x\in\overline{G_i}$ for all $i\le n$. If $x\not\in \partial G_n$, then $x\in G_n\subseteq G_i$ for all $i\le n$, consequently, $\alpha_k(x)<0$, a contradiction. Hence, $x\in\partial G_n$ and $\alpha_n(x)=0$. Conversely, if $\alpha_n(x)=0$, then $x\in \partial G_n\subseteq \overline{G_n}\subseteq\overline{G_i}$ for all $i\le n$. In consequence, $\alpha_i(x)\le 0$ for all $i<n$ and $\beta_n(x)=0$.

Now we put
\begin{gather*}
  \gamma_n(x)=\min_{i\le n}|\beta_i(x)|
\end{gather*}
and notice that
$$
\gamma_n^{-1}(0)=\bigcup_{i\le n}\partial G_i.
$$
Finally, let
\begin{gather*}
  \delta_n(x)=\left\{\begin{array}{ll}
                       0, & x\in\overline{G_n}, \\
                       \gamma_n(x), & x\in X\setminus\overline{G_n}.
                     \end{array}
  \right.
\end{gather*}
Obviously, the functions $\beta_n$, $\gamma_n$ and $\delta_n$ are continuous and $\delta_n\le \gamma_n$.

We put $U_{k,0}=F_{k,0}=\emptyset$, $k\in\mathbb N$. For all $k,n\in\mathbb N$ we define
\begin{gather*}
  U_{k,n}=\gamma_n^{-1}([0,\tfrac 1k)) \quad\mbox{and}\quad F_{k,n}=\delta_n^{-1}([\tfrac 1k,1]).
\end{gather*}

 The sets $U_{k,n}$ and $F_{k,n}$ satisfy the following conditions:
\begin{enumerate}
  \item[(A)] $U_{k,n}$ is open and $F_{k,n}$ is closed in $X$,

  \item[(B)] $\overline{U_{k+1,n}}\subseteq U_{k,n}\subseteq U_{k,n+1}$,

  \item[(C)] $F_{k,n}\subseteq F_{k+1,n}\cap F_{k,n+1}$,

  \item[(D)] $\delta_n^{-1}((0,1])=\bigcup_{k\in\mathbb N} F_{k,n}=X\setminus \bigl(\bigcup_{i=1}^n \partial G_i\cup G_n\bigr)$,

  \item[(E)] $\gamma_n^{-1}(0)=\bigcup_{i=1}^n \partial G_i=\bigcap_{k\in\mathbb N} U_{k,n}$,

  \item[(F)] $U_{k,n}\cap F_{k,n}=\emptyset$
\end{enumerate}
for all $k,n\in\mathbb N$.

Moreover, the sets $G_n$ satisfy the property
\begin{enumerate}
  \item[(G)] $(\partial G_n\setminus\partial G_{n-1})\cap\partial G_i=\emptyset$, $i<n$, $n\in\mathbb N$.
\end{enumerate}

Since the most of properties are evident, we prove only $(C)$ and $(G)$.

$(C)$. It is enough to  prove that $F_{k,n}\subseteq F_{k,n+1}$. Fix $x\in F_{k,n}$ for some $k,n\in\mathbb N$. Then $\delta_n(x)\ge \tfrac 1k$ and, in consequence, $x\not\in \overline{G_n}$. Therefore, $\alpha_n(x)>0$ and $\alpha_{n+1}(x)>0$. Hence, $\beta_n(x)>0$ and $\beta_{n+1}(x)>0$. Since $\gamma_n(x)=\delta_n(x)\ge\frac 1k$, the inequality $|\beta_i(x)|\ge \tfrac 1k$ holds for all $i\le n$. In particular, $|\beta_n(x)|=\beta_n(x)\ge\tfrac 1k$. Then $|\beta_{n+1}(x)|=\beta_{n+1}(x)\ge \beta_n(x)\ge\tfrac 1k$. Thus, $\delta_{n+1}(x)=\gamma_{n+1}(x)=\min\limits_{i\le n+1} |\beta_i(x)|\ge\tfrac 1k$ and $x\in F_{k,n+1}$.

$(G)$. Fix $x\in \partial G_n\setminus\partial G_{n-1}$. Since $x\in \partial G_n$, $x\in \overline{G_n}\setminus G_{n}$. Then $x\in\overline{G_{n-1}}$, because the sequence $(G_n)_{n\in\omega}$ is decreasing. Moreover, $x\not\in\partial G_{n-1}$ and therefore $x\in G_{n-1}$. Again, since $(G_n)_{n\in\omega}$ decreases, $x\in G_i$ for all $i<n$. Hence, $x\not \in\partial G_i$, $i<n$.

Since $X\setminus \overline{A}$ is an open subset of the $\omega$-resolvable space $X$, it is $\omega$-resolvable also. Hence, there exists a sequence $(B_k)_{k\in\omega}$ of mutually disjoint subsets of $X\setminus\overline{A}$ such that
\begin{gather*}
X\setminus\overline{A}=\bigcup_{k\in\omega} B_k
\end{gather*}
and each set $B_k$ in dense in $X\setminus \overline{A}$.

For every $x\in X$ we put
\begin{gather}
  f(x)=\left\{\begin{array}{cc}
                \frac 1n, & x\in\partial G_n\setminus \partial G_{n-1}\,\,\,\mbox{for some}\,\, n\in\mathbb N, \\
                0, & \mbox{otherwise}.
              \end{array}
  \right.
\end{gather}
Notice that $f$ is correctly defined because of property (G).

Now let
\begin{gather*}
  C_{k,n}=(U_{k,n}\setminus U_{k,n-1})\cup (B_k\cap(F_{k,n}\setminus F_{k,n-1}))
\end{gather*}
and
\begin{gather*}
  f_k(x)=\left\{\begin{array}{cc}
                  \frac 1n, & \mbox{if}\quad x\in C_{k,n} \,\,\, \mbox{for some}\,\,\,n\ge k,\\
                  0, & \mbox{otherwise}.
                \end{array}
  \right.
\end{gather*}

In order to show that $(f_k)_{k\in\omega}$ converges to $f$ pointwisely on $X$ we fix $x\in X$.

If $x\in \bigcup_{n\in\omega}\partial G_n$, then we put $N=\min\{n\in\omega: x\in\partial G_n\}$. Therefore, since $\partial G_0=\emptyset$, property (G) implies that $N\in\mathbb N$ and $x\in\partial G_N\setminus \bigcup_{i<N}\partial G_i=\partial G_{N}\setminus \partial G_{N-1}$. Hence, $f(x)=\frac 1N$. Then by (F) we conclude that $x\in\partial G_N\subseteq \bigcap_{k\in\mathbb N}U_{k,N}$ and $x\not\in\bigcap_{i<N-1}\partial G_i=\bigcap_{k\in\mathbb N}U_{k,N-1}$.
In consequence, taking into account (B) we conclude that there exists $K\in\omega$ such that $x\in U_{k,N}\setminus U_{k,N-1}$ for all $k\ge K$. Therefore, $f_k(x)=\tfrac 1N=f(x)$ for all $k\ge K$. Hence, $\lim_{k\to\infty} f_k(x)=f(x)$.

If $x\not\in \bigcup_{n\in\omega}\partial G_n$, then $f(x)=0$. Let $\varepsilon>0$ and $n\in\mathbb N$ be such that $\tfrac 1n<\varepsilon$.
Using (F) we conclude that $x\not\in\bigcup_{i\le n}\partial G_i=\bigcap_{k\in\mathbb N}U_{k,n}$. Then, taking into account (B), we obtain that there is $K_1\in\mathbb N$ such that $x\not\in U_{k,n}$ for all $k\ge K_1$. But $x\in\bigcup_{k\in\omega}B_k$ and the sets $B_k$ are disjoint, so there is $K>K_1$ such that $x\not\in \bigcup_{k\ge K} B_k$. Therefore, we have  $x\not\in \bigcup_{k\ge K} (U_{k,n}\cup B_k)$.

Assume that $k\ge K$. Consider the case $x\in U_{k,k}$. Then we choose the minimal number $m\le k$ such that $x\in U_{k,m}$. Then $m>n$ (indeed, if $m\le n$, then $x\in U_{k,m}\subseteq U_{k,n}$, a contradiction). In particular, $m>1$ and $x\in U_{k,m}\setminus U_{k,m-1}\subseteq C_{k,m}$. Therefore, $f_k(x)=\tfrac 1m<\tfrac 1n<\varepsilon$. Now we consider the case $x\not\in U_{k,k}$, then (B) implies that $x\not\in U_{k,n}$ for any $n\le k$. But $x\not\in B_k$. Hence, $x\not\in C_{k,n}$ for every $n\ge k$. Thus, $f_k(x)=0<\varepsilon$.

Now we prove that $A=UC(\mathscr F)$ for $\mathscr F=(f_n)_{n\in\omega}$.  Fix $x\in A$, $\varepsilon>0$ and let $n_0\in\mathbb N$ be such that $\tfrac{1}{n_0}<\varepsilon$. Since $x\in A$, we conclude that $x\in G_i$ for every $i$ and then $x\not\in\bigcup_{i\le n_0}\partial G_i=\bigcap_{k\in\mathbb N}U_{k,n_0}$. Then there exists $k_0\ge n_0$ with $x\not\in U_{k_0-1,n_0}$. Therefore, property (B) implies that   $x\not\in\overline{U_{k,n_0}}$ for all $k\ge k_0$.

We consider an open neighborhood
$$
U=G_{n_0}\setminus\overline{U_{k_0,n_0}}
$$
of $x$ in $X$. Take an arbitrary $u\in U$ and $k\ge k_0$. Since $u\in G_{n_0}$, we have that $u\in G_i$ and then $u\not\in\partial G_i$ for any $i\le n_0$. Therefore,
\begin{gather}
f(u)\in[0,\tfrac{1}{n_0}).
\end{gather}

Let us observe that (B) and (C) imply that $C_{k,n}\subseteq U_{k,k}\cup (B_k\cap F_{k,k})$ for all $k\ge n$. Therefore, if
$u\not\in U_{k,k}\cup (F_{k,k}\cap B_k)$, then $f_k(u)=f(u)=0$.

Let us consider the case $u\in U_{k,k}$. Then we take the minimal  $i\le k$ such that $x\in U_{k,i}$. Notice that $i>n_0$ (indeed, if $i\le n_0$, then (B) implies that $u\in U_{k,i}\subseteq U_{k,n_0}\subseteq U_{k_0,n_0}$, a contradiction). In particular, $i>1$ and then
$u\in U_{k,i}\setminus U_{k,i-1}\subseteq C_{k,i}$. Thus, $f_k(u)=\frac{1}{i}<\frac{1}{n_0}$.

In the case $u\in F_{k,k}$ we choose the minimal $j\le k$ with $u\in F_{k,j}$. Observe that $j>n_0$ (indeed, if $j\le n_0$, then $u\in F_{k,j}\subseteq F_{k,n_0}$, and so (D) implies $u\not\in G_{n_0}$, which is impossible). In particular, $j>1$ and $u\in B_k\cap(F_{k,j}\setminus F_{k,j-1})\subseteq C_{k,j}$. Therefore, $f_k(x)=\frac 1j<\frac{1}{n_0}$.

Thus, we proved that in any case $f_k(u)\in [0,\frac{1}{n_0})$.
Hence,
$$
|f(u)-f_k(u)|\le \frac{1}{n_0}<\varepsilon.
$$
Therefore, $A\subseteq UC(\mathscr F)$.

Now we prove that $UC(\mathscr F)\subseteq A$. In order to do this we fix $x\not\in A$ and show that $x\not\in UC(\mathscr F)$ in this case. Let $n_0=\max\{n\in\omega:x\in G_n\}$, $\varepsilon=\frac{1}{n_0+1}-\frac{1}{n_0+2}$, $U$ be an open neighborhood of $x$ and let $k_0\in\omega$. Notice that $x\in G_{n_0}\setminus G_{n_0+1}$.

Consider the case $x\in\overline{G_{n_0+1}}$. Then $x\in\overline{G_{n_0+1}}\setminus G_{n_0+1}=\partial G_{n_0+1}$. Therefore, (E) implies that $x\in\partial G_{n_0+1}\subseteq \bigcap_{k\in\mathbb N}U_{k,n_0+1}$. On the other hand, $x\in G_{n_0}$. So, $x\not\in\partial G_i$ for any $i\le n_0$. Using (E) we conclude that $x\not\in\bigcup_{i\le n_0}\partial G_i=\bigcap_{k\in\mathbb N}U_{k,n_0}$. Then, there exists $k_1\in\mathbb N$ such that $x\not\in U_{k_1,n_0}$. Therefore, (B) implies that $x\not\in\overline{U_{k,n_0}}$ for all $k\ge k_1$. Hence, there exists $k>\max\{k_0,n_0\}$ such that $x\in U_{k,n_0+1}\setminus\overline{U_{k,n_0}}$.

Since the set $\bigcup_{n\le n_0+1}\partial G_n$ is nowhere dense in $X$, there is a nonempty open set $V$ such that $V\subseteq U\cap (U_{k,n_0+1}\setminus\overline{U_{k,n_0}})\setminus \bigcup_{n\le n_0+1}\partial G_n$. Take $u\in V$. Then
$f(u)\in [0,\tfrac{1}{n_0+2}]$, because   $u\not \in  \bigcup_{n\le n_0+1}\partial G_n$, and
$f_k(u)=\tfrac{1}{n_0+1}$, since $u\in U_{k,n_0+1}\setminus U_{k,n_0}$. Therefore,
$$
|f(u)-f_k(u)|\ge\tfrac{1}{n_0+1}-\tfrac{1}{n_0+2}=\varepsilon.
$$

Now we assume that $x\not\in\overline{G_{n_0+1}}$. Then $x\in G_{n_0}\not\in\overline{G_{n_0+1}}$. Therefore, $x\not\in \bigcup_{i\le n_0+1}\partial G_i=\gamma_{n_0+1}^{-1}(0)$. Consequently, $\gamma_{n_0+1}(x)>0$. So, there exists  $m_0\in\mathbb N$ such that $\tfrac{1}{m_0}<\gamma_{n_0+1}(x)$. But
$x\not\in\overline{G_{n_0+1}}$. Therefore,  $\delta_{n_0+1}(x)=\gamma_{n_0+1}(x)>\tfrac{1}{m_0}$. Hence, $x\in {\rm int}\,F_{m_0,n_0+1}$.

By property (C),  there exists a number $k>\max\{k_0,m_0,n_0\}$ such that  $x\in {\rm int}F_{k,n_0+1}$. Since $x\in G_{n_0}$, $x\not \in F_{k,n_0}$. Then the set
$$
G=(U\setminus \overline{G_{n_0+1}})\cap ({\rm int} F_{k,n_0+1}\setminus F_{k,n_0})
$$
is an open neighborhood of $x$. Since $\bigcup_{n\le n_{0}+1}\partial G_n$ is nowhere dense in $X$ and $B_k$ is dense in $X\setminus \overline{A}$, there exists a point $v\in X$ such that
$$
v\in (G\setminus \bigcup_{n\le n_{0}+1}\partial G_n)\cap B_k.
$$
Then
$f_k(v)=\tfrac{1}{n_0+1}$, since $v\in C_{k,n_0+1}$, and $f(v)\in[0,\tfrac{1}{n_0+2}]$, because $v\not\in \bigcup_{n\le n_{0}+1}\partial G_n$. Hence,
$$
|f(v)-f_k(v)|\ge\tfrac{1}{n_0+1}-\tfrac{1}{n_0+2}=\varepsilon.
$$

Therefore, $A=UC(\mathscr F)$.\hfill $\Box$

  \begin{remark}
  {\rm Actually, we  use in the proof only the fact that the boundary of every open set in a topological space $X$ is a functionally closed set.
  It is find out that this is a characterization of perfectly normal spaces. Moreover, the following conditions are equivalent:
  \begin{enumerate}
    \item[(i)] $X$ is a perfectly normal space;

    \item[(ii)] every closed nowhere dense subset of $X$ is functionally closed.
  \end{enumerate}
    Evidently, (i)~$\Rightarrow$~(ii). In order to prove (ii)~$\Rightarrow$~(i) we take a closed set $F\subseteq X$. Since $\partial F=F\setminus{\rm int}\,F$ is closed and nowhere dense, there exists a continuous function $g:X\to[0,1]$ such that $\partial F=g^{-1}(0)$. Let us define $f:X\to [0,1]$ by $f(x)=g(x)$ if $x\in X\setminus {\rm int}\, F$ and $f(x)=0$ if $x\in{\rm int}\,F$. It is easy to see that $f$ is continuous and $F=f^{-1}(0)$. Therefore, $X$ is perfectly normal by Vedenisoff's theorem.}
\end{remark}

\begin{remark}
  {\rm By one of reviewers, in Theorem 1 it is sufficient to assume that  $Y$ is a non-discrete metric space.}
\end{remark}

\begin{remark}
  {\rm Any topological vector space is $\omega$-resolvable \cite{CGF}. So, the space of all continuous function $C_p([0,1])$ equipped with the topology of pointwise convergence is an example of perfectly normal $\omega$-resolvable space which is not metrizable.}
\end{remark}

\begin{remark}
  {\rm     Eric K. van Douwen proved in \cite[Theorem 5.2]{vD} that there exists a crowded countable regular space which cannot be represented as a union of two disjoint dense subsets. It is easy to see that this space is perfectly normal and  not $\omega$-resolvable.}
\end{remark}

\section{Acknowledgement}

I am very grateful to the reviewers for their careful reading of the manuscript and  valuable remarks which allowed to improve the paper.

\end{document}